\documentclass[a4paper,12pt]{amsart}
\usepackage{amsfonts}
\usepackage{amssymb}
\usepackage{ifthen}
\usepackage{amscd}
\usepackage{amsxtra}
\usepackage{graphicx}
\usepackage{color}
\nonstopmode \numberwithin{equation}{section}
\newtheorem{thm}{Theorem}%[section]
\newtheorem{lem}{Lemma}%[section]
\newtheorem{cor}{Corollary}%[section]
\newtheorem{prop}{Proposition}%[section]

\newtheorem{cl}{Claim}%[section]
\newtheorem{ca}{Case}%[section]
\newtheorem{sca}{Subcase}%[section]
\newtheorem{scl}{Subclaim}%[section]
\newtheorem{conj}{Conjecture}

\theoremstyle{definition}
\newtheorem{defn}{Definition}%[section]

\newtheorem{op}[equation]{Open Problem}
\newtheorem{ques}[equation]{Question}
\newtheorem{rem}{Remark}[section]
\newtheorem{exam}[equation]{Example}

\newcounter {own}
\def\theown {\thesection       .\arabic{own}}

\newenvironment{pf}[1][]{%
 \vskip 3mm
 \noindent
 \ifthenelse{\equal{#1}{}}%
  {{\slshape Proof. }}%
  {{\slshape #1.} }%
 }%
{\qed\bigskip}

\newcounter{alphabet}
\newcounter{tmp}
\newenvironment{Thm}[1][]{\refstepcounter{alphabet}%
\bigskip%
\noindent%
{\bf Theorem \Alph{alphabet}}%
\ifthenelse{\equal{#1}{}}{}{ (#1)}%
{\bf .} \itshape}{\vskip 8pt}

\makeatletter
\newcommand{\Ref}[1]{\@ifundefined{r@#1}{}{\setcounter{tmp}{\ref{#1}}\Alph{tmp}}}
\makeatother
\newenvironment{Lem}[1][]{\refstepcounter{alphabet}%
\bigskip%
\noindent%
{\bf Lemma \Alph{alphabet}}%
{\bf .} \itshape}{\vskip 8pt}

\newcommand{\ID}{{\mathbb D}}

\def\be{\begin{equation}}
\def\ee{\end{equation}}

\newcommand{\bee}{\begin{enumerate}}
\newcommand{\eee}{\end{enumerate}}

\newcommand{\blem}{\begin{lem}}
\newcommand{\elem}{\end{lem}}
\newcommand{\bthm}{\begin{thm}}
\newcommand{\ethm}{\end{thm}}
\newcommand{\bcor}{\begin{cor}}
\newcommand{\ecor}{\end{cor}}
\newcommand{\beg}{\begin{exam}}
\newcommand{\eeg}{\end{exam}}
\newcommand{\begs}{\begin{examples}}
\newcommand{\eegs}{\end{examples}}
\newcommand{\bdefe}{\begin{defn}}
\newcommand{\edefe}{\end{defn}}
\newcommand{\bprob}{\begin{prob}}
\newcommand{\eprob}{\end{prob}}
\newcommand{\bques}{\begin{ques}}
\newcommand{\eques}{\end{ques}}
\newcommand{\bei}{\begin{itemize}}
\newcommand{\eei}{\end{itemize}}
\newcommand{\bcon}{\begin{conj}}
\newcommand{\econ}{\end{conj}}
\newcommand{\bop}{\begin{op}}
\newcommand{\eop}{\end{op}}

\newcommand{\bca}{\begin{ca}}
\newcommand{\eca}{\end{ca}}
\newcommand{\bsca}{\begin{sca}}
\newcommand{\esca}{\end{sca}}

\newcommand{\bcl}{\begin{cl}}
\newcommand{\ecl}{\end{cl}}

\newcommand{\bscl}{\begin{scl}}
\newcommand{\escl}{\end{scl}}

\newcommand{\bcons}{\begin{conjs}}
\newcommand{\econs}{\end{conjs}}
\newcommand{\bprop}{\begin{propo}}
\newcommand{\eprop}{\end{propo}}
\newcommand{\br}{\begin{rem}}
\newcommand{\er}{\end{rem}}
\newcommand{\brs}{\begin{rems}}
\newcommand{\ers}{\end{rems}}
\newcommand{\bo}{\begin{obser}}
\newcommand{\eo}{\end{obser}}
\newcommand{\bos}{\begin{obsers}}
\newcommand{\eos}{\end{obsers}}
\newcommand{\bpf}{\begin{pf}}
\newcommand{\epf}{\end{pf}}
\newcommand{\ba}{\begin{array}}
\newcommand{\ea}{\end{array}}
\newcommand{\beq}{\begin{eqnarray}}
\newcommand{\beqq}{\begin{eqnarray*}}
\newcommand{\eeq}{\end{eqnarray}}
\newcommand{\eeqq}{\end{eqnarray*}}

\newcommand{\ds}{\displaystyle}

\newcounter{minutes}
\divide\time by 60
\newcounter{hours}
\multiply\time by 60 \addtocounter{minutes}{-\time}
\begin{document}

\bibliographystyle{amsplain}
\title[]{On characterizations of Bloch-type, Hardy-type and Lipschitz-type spaces}

%%%%%%%% BEGIN TIMESTAMP
\def\thefootnote{}
\footnotetext{ \texttt{\tiny File:~\jobname .tex,
          printed: \number\day-\number\month-\number\year,
          \thehours.\ifnum\theminutes<10{0}\fi\theminutes}
} \makeatletter\def\thefootnote{\@arabic\c@footnote}\makeatother
%%%%%%%% END TIMESTAMP

\author{Sh. Chen}
\address{Sh. Chen, Department of Mathematics and Computational
Science, Hengyang Normal University, Hengyang, Hunan 421008,
People's Republic of China.} \email{mathechen@126.com}

%\author{P. Muthukumar}
%\address{P. Muthukumar, Department of Mathematics, Indian Institute of Technology Madras,
%Chennai-600 036, India.} \email{}

\author{S. Ponnusamy $^\dagger $ %${}^{~\mathbf{*}}$
}
\address{S. Ponnusamy,
Indian Statistical Institute (ISI), Chennai Centre, SETS (Society
for Electronic Transactions and security), MGR Knowledge City, CIT
Campus, Taramani, Chennai 600 113, India. }
\email{samy@isichennai.res.in, samy@iitm.ac.in}

\author{A. Rasila }
\address{A. Rasila, Department of Mathematics and Systems Analysis, Aalto University, P. O. Box 11100, FI-00076 Aalto,
 Finland.} \email{antti.rasila@iki.fi}

%\author{X. Wang$^{\mathbf{*}}$
%${}^{~\mathbf{*}}$
%}
%\address{X. Wang, Department of Mathematics,
%Hunan Normal University, Changsha, Hunan 410081, People's Republic
%of China.} \email{xtwang@hunnu.edu.cn}

\keywords{Majorant, Banach space,  Lipschitz-type space,  Bloch-type growth  space.\\
$%{}^{\mathbf{*}}
^\dagger${\tt ~~%Corresponding author.
 This author is on leave from the Department of Mathematics,
Indian Institute of Technology Madras, Chennai-600 036, India}
%\\ ${}^{\mathbf{*}}$ Corresponding author
}
%\date{\today } %November 4, 10;
%File: Ch-W-S12${}_{}$equiv-mod${}_{}$submit.tex}

\begin{abstract}
In this paper,  we establish a Bloch-type growth theorem for
generalized Bloch-type spaces and discuss relationships between
Dirichlet-type spaces and Hardy-type spaces on certain classes of
complex-valued functions. Then we present some applications to
non-homogeneous Yukawa PDEs. We also consider some properties of the
Lipschitz-type spaces on certain classes of complex-valued
functions. Finally, we will study a class of composition operators	 on these spaces.

\end{abstract}

%\thanks{The research was partly supported by
%NSF of China (No. 11071063)} %and  Hunan Provincial Innovation
%%Foundation for Postgraduate (No. 125000-4113).  }

\maketitle \pagestyle{myheadings} \markboth{ Sh. Chen,
  S. Ponnusamy and A. Rasila}{Bloch-type  spaces,  Hardy-type spaces and  Lipschitz-type spaces}

\section{Introduction and main results }\label{csw-sec1}
%Let  $\mathbb{C}$ be the complex plane.
For $a\in\mathbb{C}$, let
$\ID(a,r)=\{z:\, |z-a|<r\}$. In particular, we use $\mathbb{D}_r$ to
denote the disk $\mathbb{D}(0,r)$ and  $\mathbb{D}$, the open unit
disk $\ID_1$. Let $\Omega$ be a domain in $\mathbb{C}$, with
non-empty boundary. Let $d_{\Omega}(z)$ be the Euclidean distance
from $z$ to the boundary $\partial \Omega$ of $\Omega$. In
particular, we always use $d(z)$ to denote the Euclidean distance
from $z$ to the boundary of $\mathbb{D}.$

%For
%$\delta_{0}>0$, let \be\label{eq2}
%\int_{0}^{\delta}\frac{\omega(t)}{t}\,dt\leq C\cdot\omega(\delta),\
%0<\delta<\delta_{0}, \ee and \be\label{eq3}
%\delta\int_{\delta}^{\infty}\frac{\omega(t)}{t^{2}}\,dt\leq
%C\cdot\omega(\delta),\ 0<\delta<\delta_{0}, \ee where $\omega$ is a
%majorant and $C$ is a positive constant. A majorant $\omega$ is said
%to be {\it regular} if it satisfies the conditions (\ref{eq2}) and
%(\ref{eq3}) (cf. \cite{D,D1,P,Pav1,Pav2}).

For a real $2\times2$ matrix $A$,
%$$A=\left(\begin{array}{cccc}
%\ds a_{11}\;~~ a_{12}\\[2mm]
%\ds a_{21}\;~~ a_{22}
%\end{array}\right),
%$$
we use the matrix norm $\|A\|=\sup\{|Az|:\,|z|=1\}$ and
the matrix function $l(A)=\inf\{|Az|:\,|z|=1\}$. With
$z=x+iy\in\mathbb{C}$,
% we define the usual %{\it
%complex differential operators%}
%$$\frac{\partial}{\partial
%z}=\frac{1}{2}\left(\frac{\partial}{\partial
%x}-i\frac{\partial}{\partial
%y}\right)~\mbox{and}~\frac{\partial}{\partial
%\overline{z}}=\frac{1}{2}\left(\frac{\partial}{\partial
%x}+i\frac{\partial}{\partial y}\right).
%$$
the formal derivative of the complex-valued functions $f=u+iv$ is given by
$$D_{f}=\left(\begin{array}{cccc}
\ds u_{x}\;~~ u_{y}\\[2mm]
\ds v_{x}\;~~ v_{y}
\end{array}\right),
$$
so that $\|D_{f}\|=|f_{z}|+|f_{\overline{z}}|$ and $l(D_{f})=\big| |f_{z}|-|f_{\overline{z}}|\big |.$
Throughout this paper, we denote by $ \mathcal{C}^{n}(\mathbb{D})$ the set of all $n$-times
continuously differentiable complex-valued function in
$\mathbb{D}$, where $n\in\{1,2,\ldots\}$.

\subsection*{Generalized Hardy spaces}
For $p\in(0,\infty]$, the {\it generalized Hardy space
$H^{p}_{g}(\mathbb{D})$} consists of all those functions
$f:\mathbb{D}\rightarrow\mathbb{C}$ such that $f$ is measurable,
$M_{p}(r,f)$ exists for all $r\in(0,1)$ and  $ \|f\|_{p}<\infty$,
where
$$\|f\|_{p}=
\begin{cases}
\displaystyle\sup_{0<r<1}M_{p}(r,f)
& \mbox{if } p\in(0,\infty),\\
\displaystyle\sup_{z\in\mathbb{D}}|f(z)| &\mbox{if } p=\infty,
\end{cases}
~\mbox{ and }~
M_{p}^{p}(r,f)=\frac{1}{2\pi}\int_{0}^{2\pi}|f(re^{i\theta})|^{p}\,d\theta.
$$
%It is not difficult to see that
The classical {\it Hardy space $H^{p}(\mathbb{D})$} consisting of analytic functions
in $\ID$ is a subspace of $H^{p}_{g}(\mathbb{D})$.

\subsection*{Generalized Bloch-type spaces}
A continuous increasing function $\omega:\, [0,\infty)\rightarrow
[0,\infty)$ with $\omega(0)=0$ is called a {\it majorant} if
$\omega(t)/t$ is non-increasing for $t>0$ (cf.
\cite{Dy1,Dy2,P,Pav1,Pav2}). Given a subset $\Omega$ of
$\mathbb{C}$, a function $f:\, \Omega\rightarrow \mathbb{C}$ is said
to belong to the {\it Lipschitz space $L_{\omega}(\Omega)$} if there
is a positive constant $C$ such that
\be\label{eq1}
|f(z)-f(w)|\leq C\omega(|z-w|) ~\mbox{ for all $z,w\in\Omega.$}
\ee

\begin{defn}\label{def-1}
For $p\in(0,\infty]$, $\alpha>0$, $\beta\in\mathbb{R}$ and a majorant $\omega$, %the
%logarithmic $\alpha$-Bloch space
we use
$\mathcal{L}_{p,\omega}\mathcal{B}^{\beta}_{\alpha}(\mathbb{D})$ to
denote the {\it generalized Bloch-type space} of all  functions $f\in \mathcal{C}^{1}(\mathbb{D})$  with
$\|f\|_{\mathcal{L}_{p,\omega}\mathcal{B}^{\beta}_{\alpha}(\mathbb{D})}<\infty$,
where
%$$\|f\|_{\mathcal{L}_{\omega}\mathcal{B}^{\beta}_{\alpha}(\mathbb{D})}=|f(0)|+\sup_{z\in\mathbb{D}}
%\left\{
%M_{p}(|z|,\|D_{f}\|)\omega\Big(d^{\alpha}(z)\Big(\log\frac{e}{d(z)}\Big)^{\beta}\Big)\right\}.$$
$$\|f\|_{\mathcal{L}_{p,\omega}\mathcal{B}^{\beta}_{\alpha}(\mathbb{D})}=
\begin{cases}
\displaystyle|f(0)|+\sup_{z\in\mathbb{D}} \left\{
M_{p}(|z|,\|D_{f}\|)\omega\Big(d^{\alpha}(z)\Big(\log\frac{e}{d(z)}\Big)^{\beta}\Big)\right\}
& \mbox{if } p\in(0,\infty),\\
\displaystyle|f(0)|+\sup_{z\in\mathbb{D}} \left\{
\|D_{f}(z)\|\omega\Big(d^{\alpha}(z)\Big(\log\frac{e}{d(z)}\Big)^{\beta}\Big)\right\}
&\mbox{if } p=\infty.
\end{cases}
$$
\end{defn}

It can be easily seen  that $\mathcal{L}_{p,\omega}\mathcal{B}^{\beta}_{\alpha}(\mathbb{D})$ is
a Banach space for $p\geq1$. Moreover, we have the following:

%\item{{{\rm(1)}}}~If we take $\omega(t)=t$, then we call $\mathcal{L}_{\omega}\mathcal{B}^{\beta}_{\alpha}(\mathbb{D})$
%as a {\it  $\beta$-logarithmic $\alpha$-Bloch space}.
\begin{enumerate}
\item[{\rm (1)}] If  $\beta=0$, then
$\mathcal{L}_{\infty,\omega}\mathcal{B}^{0}_{\alpha}(\mathbb{D})$ is
called the {\it $\omega$-$\alpha$-Bloch space}.

\item[{\rm (2)}] If we take $\alpha=1$, then
$\mathcal{L}_{\infty,\omega}\mathcal{B}^{\beta}_{1}(\mathbb{D})$ is
called  the {\it  logarithmic $\omega$-Bloch space}.

\item[{\rm (3)}] If we take $\omega(t)=t$ and $\beta=0$,
then
$\mathcal{L}_{\infty,\omega}\mathcal{B}^{0}_{\alpha}(\mathbb{D})$ is
called  the  {\it generalized $\alpha$-Bloch space} (cf. \cite{LW,
RU,Z1,Z2}).

\item[{\rm (4)}] If we take $\omega(t)=t$ and $\alpha=1$, then
$\mathcal{L}_{\infty,\omega}\mathcal{B}^{\beta}_{1}(\mathbb{D})$ is
called  the {\it generalized logarithmic Bloch space} (cf.
\cite{CPW1,Dyak,GPPJ,Pav1,Pe,Z1}).
\end{enumerate}
Let $\mathcal{A}(\mathbb{D})$ be the set of  all analytic functions
defined in $\mathbb{D}$. Then
$\mathcal{L}_{\infty,\omega}\mathcal{B}^{0}_{\alpha}(\mathbb{D})\cap\mathcal{A}(\mathbb{D})$
(resp.
$\mathcal{L}_{\infty,\omega}\mathcal{B}^{\beta}_{1}(\mathbb{D})\cap\mathcal{A}(\mathbb{D})$)
is  the {\it $\alpha$-Bloch space} (resp. {\it  logarithmic Bloch
space}), where $\omega(t)=t$.

A classical result of Hardy and Littlewood asserts that if
$p\in(0,\infty]$, $\alpha\in(1,\infty)$ and $f$ is an analytic
function in $\mathbb{D}$, then (cf. \cite{Du1,HL1,HL2})
$$ M_{p}(r,f')=O \left(\Big(\frac{1}{1-r}\Big)^{\alpha} \right ) ~\mbox{ as $r\rightarrow1$}
$$
if and only if
$$M_{p}(r,f)=O \left (\Big(\log\frac{1}{1-r}\Big )^{\alpha-1}\right) ~\mbox{ as $r\rightarrow1$}.
$$
%Indeed, the above result of Hardy and Littlewood provides a close
%relationship between the integral means of analytic functions and
%those of their derivatives (cf. \cite{Du1,HL1,HL2}).
In \cite{GPP}, Girela, Pavlovi\'c and Pel\'{a}ez  refined the above result for the
case $\alpha=1$ as follows. For related investigations in this topic, we refer to \cite{CPW2,CP,CRW,GP}.

%In \cite{GP}, Girela and Pel\'{a}ez obtained the following result.
%
\begin{Thm} $($\cite[Theorem 1.1]{GPP}$)$
\label{ThmB} Let $p\in(2,\infty)$. For $r\in(0,1)$, if  $f$ is
analytic in $\mathbb{D}$ such that
$$M_{p}(r,f')=O \left ( \frac{1}{1-r} \right  )  ~\mbox{ as $r\rightarrow1$},
$$
then
$$ M_{p}(r,f)=O \left (\Big(\log\frac{1}{1-r}\Big)^{\frac{1}{2}} \right ) ~\mbox{ as $r\rightarrow 1$}.
$$
\end{Thm}

   %In \cite{CM}
%(see also \cite[$\mbox{P}_{99}$]{KO}), Clunie and MacGregor obtained
%a similar deeper estimate under a little strong assumption, which is
%as follows.

%\begin{Thm}\label{Thm-C}
%Let $u$ be a real harmonic function in $\mathbb{D}$. For
%$r\in(0,1)$, if
%$$|\nabla u(z)|=O \left (\Big(\frac{1}{1-r}\Big)\right  )  ~\mbox{ as
%$r\rightarrow1$},
%$$
%then  $$ M_{p}(r,u)=O \left (\Big(\log\frac{1}{1-r}\Big)^{1/2}
%\right ) ~\mbox{ as $r\rightarrow 1$},$$ where $\nabla
%u=(u_{x},u_{y})$ denotes the gradient of $u$.
%\end{Thm}

\begin{defn}\label{def-3}
For $n\in\{1,2,\ldots\}$, we denote by
$\mathcal{HZ}_{n}(\mathbb{D})$ the class of all functions
$f\in\mathcal{C}^{n}(\mathbb{D})$ satisfying {\it Heinz's nonlinear
differential inequality} (cf. \cite{HZ}) \be\label{eq-15} |\Delta
f(z)|\leq a(z)\|D_{f}(z)\|+b(z)|f(z)|+q(z)~ \mbox{ {\rm
for}}~z\in\mathbb{D}, \ee where $a(z)$,  $b(z)$ and $q(z)$ are
real-valued nonnegative continuous functions in $\mathbb{D}$ and
$\Delta$ is the usual complex Laplacian operator
$$\Delta:=4\frac{\partial^{2}}{\partial z\partial
\overline{z}}=\frac{\partial^{2}}{\partial
x^{^{2}}}+\frac{\partial^{2}}{\partial y^{^{2}}}.
$$
\end{defn}

One of our primary goals is to establish a generalization of Theorem \Ref{ThmB}.

\begin{thm}\label{thm-4}
Let $\omega$ be a majorant, $p\in[2,\infty),$ $\alpha>0$,
$\beta\leq\alpha$ and $f\in \mathcal{HZ}_{2}(\mathbb{D})\cap
\mathcal{L}_{p,\omega}\mathcal{B}^{\beta}_{\alpha}(\mathbb{D})$
satisfying $\sup_{z\in\mathbb{D}}b(z)<\frac{4}{p}$,
$\sup_{z\in\mathbb{D}}a(z)<\infty$ and
$\sup_{z\in\mathbb{D}}q(z)<\infty$. If ${\rm Re}(\overline{f}\Delta
f)\geq0$, then
\begin{eqnarray*}
M_{p}(r,f)&\leq&\frac{1}{\left[1-\frac{pr^{2}}{4}\sup_{z\in\mathbb{D}}\big(b(z)\big)\right]}
\left
[\left(\frac{rp\|f\|_{\mathcal{L}_{p,\omega}\mathcal{B}^{\beta}_{\alpha}(\mathbb{D})}}{\omega(1)}\right)^{2}
\int_{0}^{1}\frac{(1-t)\,dt}
{d^{2\alpha}(rt)\Big(\log\frac{e}{d(rt)}\Big)^{2\beta}} \right .\\
&&+\frac{pr^{2}\|f\|_{\mathcal{L}_{p,\omega}\mathcal{B}^{\beta}_{\alpha}(\mathbb{D})}
\sup_{z\in\mathbb{D}}\big(a(z)\big)}{\omega(1)}M_{p}(r,f)\\
&&\left .\times\int_{0}^{1}\frac{(1-t)\,dt}
{d^{\alpha}(rt)\Big(\log\frac{e}{d(rt)}\Big)^{\beta}}
+|f(0)|^{2}+\frac{pr^{2}}{4}
\sup_{z\in\mathbb{D}}\big(q(z)\big)M_{p}(r,f)\right]^{\frac{1}{2}}.
\end{eqnarray*}
%where $C_{1}$ is a positive constant.
\end{thm}

We remark that for $\omega(t)=t$, $\alpha-1=\beta=0$ and $a(z)=
b(z)= q(z)\equiv0$, Theorem \ref{thm-4} coincides with Theorem
\Ref{ThmB}.

Let $\lambda:~\mathbb{D}\rightarrow[0,\infty)$ be continuous
and $f=u+iv$ belong  to $\mathcal{C}^{2}(\mathbb{D})$. The  elliptic partial differential equation
(or briefly the PDE) in the form
\be\label{eq-g2}
\Delta f(z)=\lambda(z) f(z)
\ee
is called the {\it non-homogeneous Yukawa PDE}. If $\lambda$ in
(\ref{eq-g2}) is a positive constant function, then we have the usual Yukawa PDE,
which first aroses from the work of the Japanese Nobel physicist Hideki Yukawa. He used this
equation to describe the nuclear potential of a point charge as
$e^{-\sqrt{\lambda} r}/r$ (cf. \cite{A,B,D-,SW,Ya}).

As an application of Theorem \ref{thm-4}, one obtains the following result.

\begin{cor}\label{cr-4}
Let $\omega$ be a majorant, $p\in[2,\infty),$
$\alpha>0$ and $\beta\leq\alpha$. Suppose
$f\in\mathcal{C}^{2}(\mathbb{D})$ and satisfies \eqref{eq-g2})
with $\sup_{z\in\mathbb{D}}\lambda(z)<\frac{4}{p}$. If
$f\in\mathcal{L}_{p,\omega}\mathcal{B}^{\beta}_{\alpha}(\mathbb{D}),$
then
\begin{eqnarray*}
M_{p}(r,f)&\leq&C_{\lambda}^{p}(r)
\left[|f(0)|^{2}+\bigg(\frac{rp\|f\|_{\mathcal{L}_{p,\omega}\mathcal{B}^{\beta}_{\alpha}(\mathbb{D})}}
{\omega(1)}\bigg)^{2}\int_{0}^{1}\frac{(1-t)\,dt}
{d^{2\alpha}(rt)\Big(\log\frac{e}{d(rt)}\Big)^{2\beta}}
\right]^{\frac{1}{2}},
\end{eqnarray*} where
$$C_{\lambda}^{p}(r)=\frac{1}{\left[1-\frac{pr^{2}}{4}\sup_{z\in\mathbb{D}}\big(\lambda(z)\big)\right]}.
$$
Furthermore, if $\alpha-1=\beta=0$, $\lambda(z)\equiv0$ is a
constant function and $\omega(t)=t$,  then
\be\label{eqw}
M_{p}(r,f)=O \left (\Big(\log\frac{1}{1-r}\Big)^{\frac{1}{2}} \right ) ~\mbox{ as $r\rightarrow 1$}
\ee
and the extremal function $f(z)=\sum_{n=0}^{\infty}z^{2^{n}}$ shows that the estimate of
\eqref{eqw} is sharp.
\end{cor}
\bpf
It is easy to see that if $f$ is a solution to (\ref{eq-g2}), then
$f$ satisfies Heinz's nonlinear differential inequality
(\ref{eq-15}). Then Corollary \ref{cr-4} follows from Theorem
\ref{thm-4}. The sharpness part in (\ref{eqw}) follows
from \cite[Theorem 1(b)]{GP}.
\epf

\begin{defn}\label{def5}
We use $\mathcal{D}_{\gamma,\mu}(\mathbb{D})$ to denote the
Dirichlet-type space consisting of  all $f\in\mathcal{C}^{1}(\mathbb{D})$ with the norm
$$\|f\|_{\mathcal{D}_{\gamma,\mu}}=|f(0)|+\int_{\mathbb{D}}d^{\gamma}(z)\|D_{f}(z)\|^{\mu}\,d\sigma(z)<\infty,
$$
where $\gamma>0$, $\mu>0$ and $d\sigma$ denotes the normalized area
measure in $\mathbb{D}$.
\end{defn}

It is not difficult to see that if $\omega(t)=t$, then
$\mathcal{L}_{1,\omega}\mathcal{B}^{0}_{\gamma}(\mathbb{D})\subset\mathcal{D}_{\gamma,1}(\mathbb{D}).$

\begin{prop}\label{prop1.1}
Let $f\in\mathcal{C}^{3}(\mathbb{D})\cap\mathcal{D}_{\gamma,2}(\mathbb{D})$
and ${\rm Re}\,[(\Delta f)_{z}\overline{f_{z}}+(\Delta
f)_{\overline{z}}\overline{f_{\overline{z}}}]\geq0.$ Then
$f\in\mathcal{L}_{\infty,\omega}\mathcal{B}^{0}_{1+\gamma/2}(\mathbb{D})$ with $\omega(t)=t.$
\end{prop}

\begin{thm}\label{thm-10}
Let
$f\in\mathcal{HZ}_{3}(\mathbb{D})\cap\mathcal{D}_{\gamma,2}(\mathbb{D})$
with ${\rm Re}\,(\overline{f}\Delta f)\geq0$ and ${\rm Re}\,[(\Delta
f)_{z}\overline{f_{z}}+(\Delta
f)_{\overline{z}}\overline{f_{\overline{z}}}]\geq0,$ where
$0<\gamma\leq1$, $\sup_{z\in\mathbb{D}}a(z)<\infty$,
$\sup_{z\in\mathbb{D}}b(z)<\infty$ and
$\sup_{z\in\mathbb{D}}q(z)<\infty$.
 If $a(z)+b(z)+q(z)$  is a non-zero
function, then $f\in H^{\frac{2}{\gamma}}_{g}(\mathbb{D})$.
\end{thm}

%Especially,  we have

%\begin{prop}\label{prop1.2}
%If $p\in(0,\infty)$ and
%$f\in\mathcal{HZ}_{2}(\mathbb{D})\cap\mathcal{D}_{0,2}(\mathbb{D})$
%with $a(z)+b(z)+q(z)\equiv0$ for $z\in\mathbb{D}$, then $f\in
%H^{p}_{g}(\mathbb{D})$.
%\end{prop}

The result given below is an easy consequence of Theorem \ref{thm-10}.

\begin{cor}\label{cor-x}
Let $0<\gamma\leq1$, $f\in\mathcal{C}^{2}(\mathbb{D})\cap\mathcal{D}_{\gamma,2}(\mathbb{D})$
and satisfy the PDE \eqref{eq-g2}, where  $\lambda(z)$ is a nonnegative
constant function.  Then $f\in H^{\frac{2}{\gamma}}_{g}(\mathbb{D})$.
\end{cor}

\subsection*{Bloch-type spaces and weighted Lipschitz functions}
Holland and  Walsh \cite{Ho}, and Zhao \cite{Zh} characterized
analytic Bloch spaces and $\alpha$-Bloch spaces in terms of weighted Lipschitz functions, respectively.
Extended discussions on this topic may be found from \cite{LW,Pav,Z1,Z2}.
Our next result characterizes  generalized $\alpha$-Bloch space by using a majorant.

\begin{thm}\label{thm-1}
Let $0\leq s<1$, $s\leq\alpha<s+1$ and $\omega$ be a majorant. Then
$f\in\mathcal{L}_{\infty,\omega}\mathcal{B}^{0}_{\alpha}(\mathbb{D})$
if and only if there is a constant $C_{1}>0$ such that, for all $z$
and $ w$ with $z\neq w$,
$$\frac{|f(z)-f(w)|}{|z-w|}\leq \frac{C_{1}}{\omega\big(d^{s}(z)d^{\alpha-s}(w)\big)}.
$$
\end{thm}

We remark that Theorem \ref{thm-1} is indeed a generalization of
\cite[Theorem 3]{Ho}, \cite[Theorem 2]{Pav} and  \cite[Theorem A]{LW} using a majorant.

\subsection*{Harmonic mappings, Bloch-type spaces and BMO}
Let $F$ be an analytic function from $\mathbb{B}^{n}$ into
$\mathbb{D}$, where $\mathbb{B}^{n}$ denotes the open unit ball in
$\mathbb{C}^{n}$. We say that $F$ has the {\it pull-back property}
if $f\circ F\in {\rm BMOA}(\mathbb{B}^{n})$ whenever analytic
function $f$ belongs to the Bloch space of $\mathbb{D}$ (cf.
\cite{RU}).

\begin{op}$($\cite[Problem 1]{RU}$)$\label{p-1}
Let $F$ be an analytic function from $\mathbb{B}^{n}$ into
$\mathbb{C}$. For which $\alpha$ does
\be\label{eq-10}\sup_{z,w\in\mathbb{B}^{n},z\neq
w}\frac{|F(z)-F(w)|}{|z-w|^{\alpha}}<\infty\ee  imply that $F$ has
the pull-back property?
\end{op}

%Let $F$ be an analytic function from  $\mathbb{B}^{n}$  into $\mathbb{C}$.
It is not difficult to see  that  $F$  satisfies
(\ref{eq-10}) if and only if
$$|\nabla F(z)|= O\big((1-|z|)^{\alpha-1}\big),
$$
where $\nabla F=(F_{z_{1}},\ldots,F_{z_{n}})$ denote the {\it complex gradient}.

A planar complex-valued function $f$ defined in $\mathbb{D}$ is
called a {\it harmonic mapping} in $D$ if and only if both the real
and the imaginary parts of $f$ are real harmonic in $\mathbb{D}$
(cf. \cite{Du}). We consider Problem \ref{p-1} for planar
harmonic mappings, and present a characterization on the relationship
between $\omega$-$\alpha$-Bloch space and ${\rm BMO}$ as follows.

\begin{thm}\label{thm-2}
Let $1\leq\alpha<2$, $f$ be a harmonic mapping in $\mathbb{D}$ and
$\omega$ be a majorant. Then
$f\in\mathcal{L}_{\infty,\omega}\mathcal{B}^{0}_{\alpha}(\mathbb{D})$
if and only if there is a constant $C_{2}>0$ such that for all
$r\in(0,d(z)]$,
$$\frac{1}{|\mathbb{D}(z,r)|}\int_{\mathbb{D}(z,r)}\left|f(\zeta)-
\frac{1}{|\mathbb{D}(z,r)|}\int_{\mathbb{D}(z,r)}f(\xi)dA(\xi)\right|dA(\zeta)\leq\frac{C_{2}r}{\omega(r^{\alpha})},
$$
where $dA$ denotes the Lebesgue area measure in $\mathbb{D}$ and
$|\mathbb{D}(z,r)|$ denotes the area of $\mathbb{D}(z,r)$.
\end{thm}

%In particular, if we take $\omega(t)=t$ in Theorem \ref{thm-2}, then
%we get

Theorem \ref{thm-2} gives the following result.

\begin{cor}\label{cor1}
Let $\alpha=1$  and $\omega$ be a majorant with $\omega(t)=t.$ Then
$f\in\mathcal{L}_{\infty,\omega}\mathcal{B}^{0}_{\alpha}(\mathbb{D})$
if and only if $f\in {\rm BMO}.$

%$f$ be a harmonic mapping in $\mathbb{D}$. Then the following are
%equivalent:

%\item{{{\rm(1)}}}~ There is a positive constant $C_{3}$ such that $$\|D_{f}(z)\|\leq \frac{C_{3}}{d^{\alpha}(z)};$$

%\item{{{\rm(2)}}}~ There is a positive constant $C_{4}$ such that for all
%$r\in(0,d(z)]$,
%$$\frac{1}{|\mathbb{D}(z,r)|}\int_{\mathbb{D}(z,r)}\left|f(\zeta)-
%\frac{1}{|\mathbb{D}(z,r)|}\int_{\mathbb{D}(z,r)}f(\xi)dA(\xi)\right|dA(\zeta)\leq
%C_{4}r^{2(1-\alpha)}.
\end{cor}

By Theorems \ref{thm-1} and \ref{thm-2}, we also have the following.

\begin{cor}\label{cor2}
Let $0\leq s<1$,  $1\leq\alpha<s+1$ and $f$ be a harmonic mapping in $\mathbb{D}$.
Then the following are equivalent:
\begin{enumerate}
\item[{\rm (1)}] $f\in\mathcal{L}_{\infty,\omega}\mathcal{B}^{0}_{\alpha}(\mathbb{D});$

\item[{\rm (2)}] There exists a constant $C_4>0$ such that for all $z, w\in\mathbb{D}$ with $z\neq w$,
$$\frac{|f(z)-f(w)|}{|z-w|}\leq \frac{C_{4}}{\omega\big(d^{s}(z)d^{\alpha-s}(w)\big)};
$$

\item[{\rm (3)}] There exists a constant $C_5>0$ such that for all  $r\in(0,d(z)]$,
$$\frac{1}{|\mathbb{D}(z,r)|}\int_{\mathbb{D}(z,r)}\left|f(\zeta)-
\frac{1}{|\mathbb{D}(z,r)|}\int_{\mathbb{D}(z,r)}f(\xi)\,dA(\xi)\right|dA(\zeta)\leq\frac{C_{5}r}{\omega(r^{\alpha})}.
$$
\end{enumerate}
\end{cor}

\begin{defn}\label{def-2}
The {\it little Bloch-type space}
$\mathcal{L}_{\infty,\omega}^{0}\mathcal{B}^{\beta}_{\alpha}(\mathbb{D})$
consists of all functions $f\in\mathcal{L}_{\infty,\omega}\mathcal{B}^{\beta}_{\alpha}(\mathbb{D})$
such that
$$\lim_{|z|\rightarrow1-}\left\{
\|D_{f}(z)\|\omega\Big(d^{\alpha}(z)\Big(\log\frac{e}{d(z)}\Big)^{\beta}\Big)\right\}=0.
$$
\end{defn}

Our next result provides a characterization for the little Bloch-type space
$\mathcal{L}_{\infty,\omega}^{0}\mathcal{B}^{0}_{\alpha}(\mathbb{D}).$

\begin{thm}\label{thm-3}
Let $0\leq s<1$, $s\leq\alpha<s+1$ and $\omega$ be a majorant. Then
$f\in\mathcal{L}_{\infty,\omega}^{0}\mathcal{B}^{0}_{\alpha}(\mathbb{D})$
if and only if
\be\label{eq12}\lim_{|z|\rightarrow1-}\sup_{ w\in\mathbb{D},z\neq
w}\left\{\frac{|f(z)-f(w)|\omega\big(d^{s}(z)d^{\alpha-s}(w)\big)}{|z-w|}\right\}=0.
\ee
\end{thm}

\subsection*{Composition operators}
 If $\omega(t)=t$, then we denote
$\mathcal{L}\mathcal{B}^{\beta}_{\alpha}(\mathbb{D})=\mathcal{A}(\mathbb{D})
\cap\mathcal{L}_{\infty,\omega}\mathcal{B}^{\beta}_{\alpha}(\mathbb{D}).$
Given an analytic self mapping  $\phi$ of the unit disk $\mathbb{D}$,
the composition operator $C_{\phi}:\,\mathcal{A}(\mathbb{D})\rightarrow
\mathcal{A}(\mathbb{D})$ is defined by
$$C_{\phi}(f)=f\circ\phi,
$$
where $f\in \mathcal{A}(\mathbb{D})$ (cf. \cite{AD,Pav1,Pe,Shap,Z1}).

\begin{thm}\label{thm-8}
Let $\alpha>0$, $\beta\leq\alpha$ and
$\phi:\,\mathbb{D}\rightarrow\mathbb{D}$  be an analytic function.
Then the following are equivalent:
\begin{enumerate}
\item[{\rm (1)}] $C_{\phi}:\,\mathcal{L}\mathcal{B}^{\beta}_{\alpha}(\mathbb{D})\rightarrow H^{2}(\mathbb{D})$
is a bounded operator;

\item[{\rm (2)}] $\displaystyle \frac{1}{2\pi}\int_{0}^{2\pi}\int_{0}^{1}
\frac{|\phi'(re^{i\theta})|^{2}}{d^{2\alpha}(\phi(re^{i\theta}))}
\left(\log\frac{e}{d(\phi(re^{i\theta}))}\right)^{-2\beta}(1-r)\,dr\,d\theta<\infty.
$
\end{enumerate}
\end{thm}

The proofs of Theorems \ref{thm-4} and \ref{thm-10}  will be
presented in Section \ref{csw-sec2}, and the proofs of Theorems
\ref{thm-1}, \ref{thm-2} and \ref{thm-3}   will be given in Section
\ref{csw-sec3}. Theorem \ref{thm-8} will be proved in the last
section.

\section{Bloch-type growth spaces and applications to PDEs}\label{csw-sec2}

Green's theorem (cf. \cite{Pav3}) states that if $g\in
\mathcal{C}^{2}(\mathbb{D})$,  then for $r\in (0, 1)$,
\be\label{eq1.2x}
\frac{1}{2\pi}\int_{0}^{2\pi}g(re^{i\theta})\,d\theta=g(0)+
\frac{1}{2}\int_{\mathbb{D}_{r}}\Delta
g(z)\log\frac{r}{|z|}\,d\sigma(z).
\ee

\begin{lem}\label{lem-3}
Let $f\in\mathcal{C}^{2}(\mathbb{D})$ such that ${\rm
Re}\,(\overline{f}\Delta f)\geq0$. Then for $p\in[2,\infty)$,
$M_{p}^{p}(r,f)$ is an increasing function of $r$, $r\in(0,1).$
\end{lem}
\bpf First we deal with the case $p\in[2,4)$. In this case,
%\subsection*{Proof of Lemma \ref{lem-3}}
%\bca Let $p\in[2,4)$.\label{case-1} \eca
for $n\in\{1,2,\ldots\}$, we let $F_{n}^{p}=\left(|f|^{2}+\frac{1}{n}\right)^{\frac{p}{2}}$.
Then, by elementary calculations, we have
\begin{eqnarray*}
\Delta(F_{n}^{p})&=&4\frac{\partial^{2}}{\partial
z\partial\overline{z}}(F_{n}^{p})\\
&=&p(p-2)\left(|f|^{2}+\frac{1}{n}\right)^{\frac{p}{2}-2}|f\overline{f_{z}}
+f_{\overline{z}}\overline{f}|^{2}\\
&&+2p\left(|f|^{2}+\frac{1}{n}\right)^{\frac{p}{2}-1}\big(|f_{z}|^{2}+|f_{\overline{z}}|^{2}\big)
+p\left(|f|^{2}+\frac{1}{n}\right)^{\frac{p}{2}-1}{\rm Re}\,(\overline{f}\Delta f).
\end{eqnarray*}
Let $\tau_{n}=\Delta(F_{n}^{p})$ and
$$\tau=p(p-2)|f|^{p-2}\|D_{f}\|^{2}+2p\left(|f|^{2}+1\right)^{\frac{p}{2}-1}\big(|f_{z}|^{2}+|f_{\overline{z}}|^{2}\big)
+p\left(|f|^{2}+1\right)^{\frac{p}{2}-1}{\rm Re}\,(\overline{f}\Delta f).
$$
For $r\in(0,1)$, it is not difficult to see that $\tau_{n}$
and $\tau$ are integrable in $\mathbb{D}_{r}$, and
$\tau_{n}\leq\tau$.

By (\ref{eq1.2x}) and Lebesgue's dominated convergence theorem, we
conclude that
\begin{eqnarray*}
\lim_{n\rightarrow\infty}r\frac{d}{dr}M_{p}^{p}(r,
F_{n})&=&\frac{1}{2}\lim_{n\rightarrow\infty}\int_{\mathbb{D}_{r}}\tau_{n}(z)\,d\sigma(z)\\
&=&\frac{1}{2}\int_{\mathbb{D}_{r}}\lim_{n\rightarrow\infty}\tau_{n}(z)\,d\sigma(z)\\
&=&\frac{1}{2}\int_{\mathbb{D}_{r}}\Big[p(p-2)|f(z)|^{p-4}|f(z)\overline{f_{z}(z)}+
f_{\overline{z}}(z)\overline{f(z)}|^{2}\\
&&+2p|f(z)|^{p-2}\big(|f_{z}(z)|^{2}+|f_{\overline{z}}(z)|^{2}\big)\\
&&+p|f(z)|^{p-2}{\rm Re}\big(\overline{f(z)}\Delta
f(z)\big)\Big]\,d\sigma(z)\\
&=&r\frac{d}{dr}M_{p}^{p}(r,f),
\end{eqnarray*}
which implies that $M_{p}^{p}(r,f)$ is increasing with respect to $r$ in $(0,1)$.

%\bca Let $p\in[4,\infty)$.\label{case-2} \eca
Next we consider the case $p\in[4,\infty)$. Since
\begin{eqnarray*}
\Delta(|f|^{p})&=&p(p-2)|f|^{p-4}|f\overline{f_{z}}
+f_{\overline{z}}\overline{f}|^{2}\\
&&+2p|f|^{p-2}\big(|f_{z}|^{2}+|f_{\overline{z}}|^{2}\big)
+p|f|^{p-2}{\rm Re}\,(\overline{f}\Delta f)\geq0,
\end{eqnarray*}
we see that $|f|^{p}$ is subharmonic in $\mathbb{D}$. Hence
$M_{p}^{p}(r,f)$ is also increasing with respect to $r\in(0,1)$, and the proof is complete.
\epf%\qed

\begin{lem}\label{lem-4}
Let $f\in\mathcal{C}^{2}(\mathbb{D})$ with ${\rm
Re}\,(\overline{f}\Delta f)\geq0.$ Then for $p\in[2,\infty)$,
$$\int_{\mathbb{D}_{r}}|f(z)|^{p}\log\frac{r}{|z|}\,d\sigma(z)\leq\frac{r^{2}}{2}M_{p}^{p}(r,f).
$$
\end{lem}
%\subsection*{Proof of Lemma \ref{lem-4}}
\bpf By Lemma \ref{lem-3}, we see that
\begin{eqnarray*}
\int_{\mathbb{D}_{r}}|f(z)|^{p}\log\frac{r}{|z|}\,d\sigma(z)&=&\frac{1}{\pi}\int_{0}^{2\pi}\int_{0}^{r}
|f(\rho e^{i\theta})|^{p}\rho\log\frac{r}{\rho}\,d\rho \,d\theta\\
&=&2\int_{0}^{r}M_{p}^{p}(\rho,f)\rho\log\frac{r}{\rho}\,d\rho\\
&\leq&2M_{p}^{p}(r,f)\int_{0}^{r}\rho\log\frac{r}{\rho}\,d\rho\\
&=&\frac{r^{2}}{2}M_{p}^{p}(r,f).
\end{eqnarray*}
The proof of the lemma is complete.
\epf %\qed

The following lemma easily follows from elementary computations
and the monotonicity of the function $\omega(t)/t$.

\begin{lem}\label{lem-5}
Suppose that $\alpha>0$, $\beta\leq\alpha$ and $\omega$ is a majorant.
For $r\in(0,1)$, let
$$\eta(r)=d^{\alpha}(r)\left(\log\frac{e}{d(r)}\right)^{\beta}.
$$
Then $\eta(r)$  and $\eta(r)/\omega(\eta(r))$ are decreasing in
$(0,1)$.
\end{lem}

\subsection*{Proof of Theorem \ref{thm-4}}
By H\"older's inequality, we have

 \beq\label{eq-cpr1}
\frac{1}{2\pi}\int_{0}^{2\pi}|f(re^{i\theta})|^{p-1}\|D_{f}(re^{i\theta})\|\,d\theta
&\leq&\left(\frac{1}{2\pi}\int_{0}^{2\pi}|f(re^{i\theta})|^{p}\,d\theta\right)^{\frac{p-1}{p}}\\
\nonumber &&\times\left(\frac{1}{2\pi}\int_{0}^{2\pi}
\|D_{f}(re^{i\theta})\|^{p}\,d\theta\right)^{\frac{1}{p}}\\
\nonumber &=&M_{p}^{p-1}(r,f)M_{p}(r,\|D_{f}\|),
\eeq
\beq\label{eq-cpr2}
\frac{1}{2\pi}\int_{0}^{2\pi}|f(re^{i\theta})|^{p-2}\|D_{f}(re^{i\theta})\|^{2}\,d\theta
&\leq&
\left(\frac{1}{2\pi}\int_{0}^{2\pi}|f(re^{i\theta})|^{p}\,d\theta\right)^{\frac{p-2}{p}}\\
\nonumber&&\times\left(\frac{1}{2\pi}\int_{0}^{2\pi}
\|D_{f}(re^{i\theta})\|^{p}\,d\theta\right)^{\frac{2}{p}}\\
\nonumber &=&M_{p}^{p-2}(r,f)M_{p}^{2}(r,\|D_{f}\|). \eeq
and
\beq\label{eq-cpr3}
\frac{1}{2\pi}\int_{0}^{2\pi}|f(re^{i\theta})|^{p-1}\,d\theta
&\leq&\left(\frac{1}{2\pi}\int_{0}^{2\pi}|f(re^{i\theta})|^{p}\,d\theta\right)^{\frac{p-1}{p}}\\
\nonumber &&\times\left(\frac{1}{2\pi}\int_{0}^{2\pi}
\,d\theta\right)^{\frac{1}{p}}\\
\nonumber &=&M_{p}^{p-1}(r,f). \eeq

 By (\ref{eq-15}),
(\ref{eq1.2x}), (\ref{eq-cpr1}), (\ref{eq-cpr2}), (\ref{eq-cpr3}),
Lemmas \ref{lem-4} and \ref{lem-5}, and Lebesgue's dominated
convergence theorem, we see that
\begin{eqnarray*}
M_{p}^{p}(r,f)&=&|f(0)|^{p}+\frac{1}{2}\int_{\mathbb{D}_{r}}\Delta\big(|f(z)|^{p}\big)\log\frac{r}{|z|}\,d\sigma(z)\\
&=&|f(0)|^{p}+\frac{1}{2}\int_{\mathbb{D}_{r}}\Big[p(p-2)|f(z)|^{p-4}|f(z)\overline{f_{z}(z)}
+f_{\overline{z}}(z)\overline{f(z)}|^{2}\\
&&+2p|f(z)|^{p-2}\big(|f_{z}(z)|^{2}+|f_{\overline{z}}(z)|^{2}\big)\\
&&+p|f(z)|^{p-2}{\rm Re}\big(\overline{f(z)}\Delta f(z)\big)
\Big]\log\frac{r}{|z|}\,d\sigma(z)\\
&\leq&|f(0)|^{p}+\frac{1}{2}\int_{\mathbb{D}_{r}}\Big(p^{2}|f(z)|^{p-2}\|D_{f}(z)\|^{2}\\
&&+p|f(z)|^{p-1}|\Delta f(z)|\Big)\log\frac{r}{|z|}\,d\sigma(z)
%\\
\end{eqnarray*}
\begin{eqnarray*}
&\leq&|f(0)|^{p}+\frac{p}{2}\int_{\mathbb{D}_{r}}\Big(p|f(z)|^{p-2}\|D_{f}(z)\|^{2}+b(z)|f(z)|^{p}\\
&&+a(z)|f(z)|^{p-1}\|D_{f}(z)\|+q(z)|f(z)|^{p-1}\Big)\log\frac{r}{|z|}\,d\sigma(z)\\
&=&|f(0)|^{p}+p^{2}\int_{0}^{r}\left(\frac{1}{2\pi}\int_{0}^{2\pi}|f(\rho
e^{i\theta})|^{p-2}\|D_{f}(\rho
e^{i\theta})\|^{2}d\theta\right)\rho\log\frac{r}{\rho}\,d\rho\\
&&+p\sup_{z\in\mathbb{D}}\big(a(z)\big)\int_{0}^{r}\left(\frac{1}{2\pi}\int_{0}^{2\pi}|f(\rho
e^{i\theta})|^{p-1}\|D_{f}(\rho
e^{i\theta})\|d\theta\right)\rho\log\frac{r}{\rho}\,d\rho\\
&&+p\sup_{z\in\mathbb{D}}\big(b(z)\big)\int_{0}^{r}\left(\frac{1}{2\pi}\int_{0}^{2\pi}|f(\rho
e^{i\theta})|^{p}d\theta\right)\rho\log\frac{r}{\rho}\,d\rho
\\
&&+p\sup_{z\in\mathbb{D}}\big(q(z)\big)\int_{0}^{r}\left(\frac{1}{2\pi}\int_{0}^{2\pi}|f(\rho
e^{i\theta})|^{p-1}d\theta\right)\rho\log\frac{r}{\rho}\,d\rho\\
&\leq&|f(0)|^{p}+p^{2}\int_{0}^{r}M_{p}^{p-2}(\rho,f)M_{p}^{2}(\rho,\|D_{f}\|)\rho\log\frac{r}{\rho}\,d\rho\\
&&+p\sup_{z\in\mathbb{D}}\big(a(z)\big)\int_{0}^{r}M_{p}^{p-1}(\rho,f)M_{p}(\rho,\|D_{f}\|)\rho\log\frac{r}{\rho}\,d\rho\\
&&+\frac{pr^{2}}{4}\sup_{z\in\mathbb{D}}\big(b(z)\big)M_{p}^{p}(r,f)\\
&&+p\sup_{z\in\mathbb{D}}\big(q(z)\big)\int_{0}^{r}M_{p}^{p-1}(\rho,f)\rho\log\frac{r}{\rho}\,d\rho\\
&\leq&|f(0)|^{p}+p^{2}\int_{0}^{r}M_{p}^{p-2}(\rho,f)M_{p}^{2}(\rho,\|D_{f}\|)\rho\log\frac{r}{\rho}\,d\rho\\
&&+p\sup_{z\in\mathbb{D}}\big(a(z)\big)\int_{0}^{r}M_{p}^{p-1}(\rho,f)M_{p}(\rho,\|D_{f}\|)\rho\log\frac{r}{\rho}\,d\rho\\
&&+\frac{pr^{2}}{4}\sup_{z\in\mathbb{D}}\big(b(z)\big)M_{p}^{p}(r,f)+\frac{pr^{2}}{4}
\sup_{z\in\mathbb{D}}\big(q(z)\big)M_{p}^{p-1}(r,f)
\end{eqnarray*}
which gives
\begin{eqnarray*}
C_{b}^{p}(r)M_{p}^{2}(r,f) &=&
\left[1-\frac{pr^{2}}{4}\sup_{z\in\mathbb{D}}\big(b(z)\big)\right]M_{p}^{2}(r,f)
\\
&\leq&|f(0)|^{2}+p^{2}\int_{0}^{r}M_{p}^{2}(\rho,\|D_{f}\|)\rho\log\frac{r}{\rho}\,d\rho\\
&&+p\sup_{z\in\mathbb{D}}\big(a(z)\big)\int_{0}^{r}M_{p}(\rho,f)M_{p}(\rho,\|D_{f}\|)\rho\log\frac{r}{\rho}\,d\rho
\\
&&+\frac{pr^{2}}{4}
\sup_{z\in\mathbb{D}}\big(q(z)\big)M_{p}(r,f)
%\\
\end{eqnarray*}
\begin{eqnarray*}
&=&|f(0)|^{2}+p^{2}\int_{0}^{r}M_{p}^{2}(\rho,\|D_{f}\|)(r-\rho)\,d\rho\\
&&+p\sup_{z\in\mathbb{D}}\big(a(z)\big)M_{p}(r,f)\int_{0}^{r}M_{p}(\rho,\|D_{f}\|)(r-\rho)dt\\
&&+\frac{pr^{2}}{4} \sup_{z\in\mathbb{D}}\big(q(z)\big)M_{p}(r,f)\\
&=&|f(0)|^{2}+(rp)^{2}\int_{0}^{1}M_{p}^{2}(tr,\|D_{f}\|)(1-t)\,dt\\
&&+pr^{2}\sup_{z\in\mathbb{D}}\big(a(z)\big)M_{p}(r,f)\int_{0}^{1}M_{p}(tr,\|D_{f}\|)(1-t)\,dt\\
&&+\frac{pr^{2}}{4} \sup_{z\in\mathbb{D}}\big(q(z)\big)M_{p}(r,f)\\
&\leq&|f(0)|^{2}+\big(rp\|f\|_{\mathcal{L}_{p,\omega}\mathcal{B}^{\beta}_{\alpha}(\mathbb{D})}\big)^{2}
\int_{0}^{1}\frac{d^{2\alpha}(rt)\Big(\log\frac{e}{d(rt)}\Big)^{2\beta}}
{\omega^{2}\left(d^{\alpha}(rt)\Big(\log\frac{e}{d(rt)}\Big)^{\beta}\right)}\\
&&\times\frac{(1-t)\,dt}{d^{2\alpha}(rt)\Big(\log\frac{e}{d(rt)}\Big)^{2\beta}}
+ pr^{2}\|f\|_{\mathcal{L}_{p,\omega}\mathcal{B}^{\beta}_{\alpha}(\mathbb{D})}
\sup_{z\in\mathbb{D}}\big(a(z)\big)M_{p}(r,f)\\
&&\times\int_{0}^{1}\frac{d^{\alpha}(rt)\Big(\log\frac{e}{d(rt)}\Big)^{\beta}}
{\omega\left(d^{\alpha}(rt)\Big(\log\frac{e}{d(rt)}\Big)^{\beta}\right)}
\frac{(1-t)\,dt}{d^{\alpha}(rt)\Big(\log\frac{e}{d(rt)}\Big)^{\beta}}
\\&&+\frac{pr^{2}}{4} \sup_{z\in\mathbb{D}}\big(q(z)\big)M_{p}(r,f)\\
&\leq&|f(0)|^{2}+\left(\frac{rp\|f\|_{\mathcal{L}_{p,\omega}\mathcal{B}^{\beta}_{\alpha}(\mathbb{D})}}
{\omega(1)}\right)^{2}\int_{0}^{1}\frac{(1-t)\,dt} {d^{2\alpha}(rt)\Big(\log\frac{e}{d(rt)}\Big)^{2\beta}}
\\
&&+\frac{pr^{2}\|f\|_{\mathcal{L}_{p,\omega}\mathcal{B}^{\beta}_{\alpha}(\mathbb{D})}
\sup_{z\in\mathbb{D}}\big(a(z)\big)}{\omega(1)}M_{p}(r,f)
%\\&& \times
\int_{0}^{1}\frac{(1-t)\,dt}
{d^{\alpha}(rt)\Big(\log\frac{e}{d(rt)}\Big)^{\beta}}\\
&&+\frac{pr^{2}}{4} \sup_{z\in\mathbb{D}}\big(q(z)\big)M_{p}(r,f),
\end{eqnarray*}
where
$$C_{b}^{p}(r)=1-\frac{pr^{2}}{4}\sup_{z\in\mathbb{D}}\big(b(z)\big).
$$
%Hence
%\begin{eqnarray*}
%M_{p}(r,f)&\leq&\frac{1}{C_{b}(p)}
%\bigg[\left(\frac{rp\|f\|_{\mathcal{L}_{p,\omega}\mathcal{B}^{\beta}_{\alpha}(\mathbb{D})}}
%{\omega(1)}\right)^{2}\int_{0}^{1}\frac{(1-t)\,dt}
%{d^{2\alpha}(rt)\Big(\log\frac{e}{d(rt)}\Big)^{2\beta}}\\
%&&+\frac{pr^{2}\|f\|_{\mathcal{L}_{p,\omega}\mathcal{B}^{\beta}_{\alpha}(\mathbb{D})}
%\sup_{z\in\mathbb{D}}\big(a(z)\big)}{\omega(1)}M_{p}(r,f)\\
%&&\times\int_{0}^{1}\frac{(1-t)\,dt}
%{d^{\alpha}(rt)\Big(\log\frac{e}{d(rt)}\Big)^{\beta}}+|f(0)|^{2}\bigg]^{\frac{1}{2}}.
%\end{eqnarray*}
%The proof of this theorem is complete.
The desired conclusion follows.
\qed

%\subsection*{Proof of Corollary \ref{cr-4}}
%It is easy to see that if $f$ is a solution to (\ref{eq-g2}), then
%$f$ satisfies Heinz's nonlinear differential inequality
%(\ref{eq-15}). Then Corollary \ref{cr-4} follows from Theorem
%\ref{thm-4}. The sharpness part in (\ref{eqw}) follows
%from \cite[Theorem 1(b)]{GP}. \qed

\begin{lem}\label{lem-9}
Let $f\in\mathcal{C}^{3}(\mathbb{D})$ with ${\rm Re}\,[(\Delta
f)_{z}\overline{f_{z}}+(\Delta
f)_{\overline{z}}\overline{f_{\overline{z}}}]\geq0.$ Then
$F=|f_{z}|^{2}+|f_{\overline{z}}|^{2}$ is subharmonic in
$\mathbb{D}$.
\end{lem}
%\subsection*{Proof of Lemma \ref{lem-9}}
\bpf Since
$F_z=f_{zz}\overline{f_{z}}+f_{z}\overline{f_{z\overline{z}}}+f_{\overline{z}z}\overline{f_{\overline{z}}}
+f_{\overline{z}}\overline{f_{\overline{z}\overline{z}}},
$
we see that
$$\Delta F=4\frac{\partial^{2}F}{\partial z\partial
\overline{z}}=4(|f_{zz}|^{2}+|f_{\overline{z}\overline{z}}|^{2})+\frac{1}{2}|\Delta
f|^{2}+2{\rm Re}\,[(\Delta f)_{z}\overline{f_{z}}+(\Delta
f)_{\overline{z}}\overline{f_{\overline{z}}}]\geq0.
$$
Then $F$ is subharmonic in $\mathbb{D}$.
%The proof of this lemma is complete.
\epf %\qed

\subsection*{Proof of Proposition \ref{prop1.1}}

By Lemma \ref{lem-9}, we know that
$F=|f_{z}|^{2}+|f_{\overline{z}}|^{2}$ is subharmonic in
$\mathbb{D}$. Then for $r\in[0,d(z))$, we have
$$F(z)\leq\frac{1}{2\pi}\int_{0}^{2\pi} F(z+re^{i\theta})\,d\theta.
$$
Integration leads to
\begin{eqnarray*}
\frac{ d^{2}(z)F(z)}{4}
&\leq&\int_{0}^{2\pi}\int_{0}^{\frac{d(z)}{2}}r|F
(z+re^{i\theta})|\frac{drd\theta}{\pi}\\
&=&\int_{\mathbb{D}(z,\frac{d(z)}{2})}
F(\zeta)\,d\sigma(\zeta)\\
%\end{eqnarray*}\begin{eqnarray*}
&\leq&2^{\gamma}d^{-\gamma}(z)\int_{\mathbb{D}(z,\frac{d(z)}{2})}d^{\gamma}(\zeta)
F(\zeta)\,d\sigma(\zeta)\\
&\leq&2^{\gamma} \|f\|_{\mathcal{D}_{\gamma,2}}d^{-\gamma}(z),
\end{eqnarray*}
which gives
\be\label{eq-u1}
\|D_{f}(z)\|\leq\sqrt{2F(z)}\leq\frac{C_{6}}{(d(z))^{1+\gamma/2}},
\ee
where $C_{6}=2^{\frac{\gamma+3}{2}} \sqrt{\|f\|_{\mathcal{D}_{\gamma,2}}}$.
Hence
$$\sup_{z\in\mathbb{D}}\left\{(d(z))^{1+\gamma/2}\|D_{f}(z)\|\right\}<\infty,
$$
which implies that
$f\in\mathcal{L}_{\infty,\omega}\mathcal{B}^{0}_{1+\gamma/2}(\mathbb{D})$,
where $\omega(t)=t.$   \qed

\vspace{6pt}

The following result is well-known.
\begin{lem}\label{Lemx}
Suppose that $a,b\in[0,\infty)$ and $q\in(0,\infty)$. Then
$$(a+b)^{q}\leq2^{\max\{q-1,0\}}(a^{q}+b^{q}).
$$
\end{lem}

\subsection*{Proof of Theorem \ref{thm-10}}
We first prove that
\be\label{eq-p1}
\int_{\mathbb{D}}d(z)\Delta(|f(z)|^{2/\gamma})\,d\sigma(z)<\infty.
\ee
 %By Lemma \ref{lem-9}, we know that
%$F=|f_{z}|^{2}+|f_{\overline{z}}|^{2}$ is subharmonic in
%$\mathbb{D}$. Then for $r\in[0,d(z))$, we have
%$$F(z)\leq\frac{1}{2\pi}\int_{0}^{2\pi}
%F(z+re^{i\theta})\,d\theta.
%$$
%Integration leads to
%\begin{eqnarray*}
%\frac{\pi d^{2}(z)F(z)}{4}
%&\leq&\int_{0}^{2\pi}\int_{0}^{\frac{d(z)}{2}}r|F
%(z+re^{i\theta})|\,dr\,d\theta\\
%&=&\pi\int_{\mathbb{D}(z,\frac{d(z)}{2})}
%F(\zeta)\,d\sigma(\zeta)\\
%\end{eqnarray*}\begin{eqnarray*}
%&\leq&\pi2^{\gamma}d^{-\gamma}(z)\int_{\mathbb{D}(z,\frac{d(z)}{2})}d^{\gamma}(\zeta)
%F(\zeta)\,d\sigma(\zeta)\\
%&\leq&\pi2^{\gamma} \|f\|_{\mathcal{D}_{\gamma}}d^{-\gamma}(z),
%\end{eqnarray*}
%which gives
%$$\|D_{f}(z)\|\leq\sqrt{2F(z)}\leq\frac{M_{1}}{(d(z))^{1+\gamma/2}},
%$$
%where $M_{1}=2^{\gamma+\frac{5}{2}} \|f\|_{\mathcal{D}_{\gamma}}$.
By (\ref{eq-u1}), we have
\begin{eqnarray*}
|f(z)|&\leq& |f(0)|+\left|\int_{[0,z]}\,df(\zeta)\right|\\
&\leq&|f(0)|+\int_{[0,z]}\|D_{f}(\zeta)\|\,|d\zeta|\\
&\leq&|f(0)|+\frac{C_{7}}{(d(z))^{\gamma/2}},
\end{eqnarray*}
where $C_{7}=\Big(2^{\frac{\gamma+5}{2}}
\sqrt{\|f\|_{\mathcal{D}_{\gamma,2}}}\Big)/\gamma$ and $[0,z]$
denotes the line segment from $0$ to $z$. Let $p=2/\gamma$. Then
Lemma \ref{Lemx} implies that for $z\in\ID$, \be\label{eq-20}
|f(z)|^{p}\leq\left[|f(0)|+\frac{C_{7}}{(d(z))^{1/p}}\right]^{p}
\leq2^{p-1}\left[|f(0)|^{p}+\frac{C_{7}^{p}}{d(z)}\right], \ee
\be\label{eq-30}
|f(z)|^{p-1}\leq\left[|f(0)|+\frac{C_{7}}{(d(z))^{1/p}}\right]^{p-1}
\leq2^{p-2}\left[|f(0)|^{p-1}+\frac{C_{7}^{p-1}}{(d(z))^{(p-1)/p}}\right]
\ee and \be\label{eq-40}
|f(z)|^{p-2}\leq\left[|f(0)|+\frac{C_{7}}{(d(z))^{1/p}}\right]^{p-2}
\leq2^{p-2}\left[|f(0)|^{p-2}+\frac{C_{7}^{p-2}}{(d(z))^{(p-2)/p}}\right].
\ee

We divide the remaining part of the proof into two cases, namely $p\in[4,\infty)$ and $p\in[2,4)$.
For the case $p\in[4,\infty)$, easy calculations give
\begin{align*}
\Delta(|f|^{p})& =4\frac{\partial^{2}}{\partial
z\partial\overline{z}}(|f|^{p})\\
&\leq p^{2}|f|^{p-2}\|D_{f}\|^{2}+p|f|^{p-1}|\Delta f|\\
&\leq
p^{2}|f|^{p-2}\|D_{f}\|^{2}+pa|f|^{p-1}\|D_{f}\|+pb|f|^{p}+pq|f|^{p-1}.
\end{align*}
Hence we infer from  (\ref{eq-20}),  (\ref{eq-30}) and (\ref{eq-40})
that for $z\in\ID$,
\beq\label{eq-50}
\nonumber d(z)\Delta(|f(z)|^{p})&\leq&p^{2}d(z)|f(z)|^{p-2}\|D_{f}(z)\|^{2}+pq(z)|f(z)|^{p-1}\\
\nonumber &&+pad(z)|f(z)|^{p-1}\|D_{f}(z)\|+pb(z)d(z)|f(z)|^{p}\\
\nonumber
&=&p^{2}(d(z))^{1-\frac{2}{p}}|f(z)|^{p-2}(d(z))^{\frac{2}{p}}\|D_{f}(z)\|^{2}\\
\nonumber
&&+p\sup_{z\in\mathbb{D}}(a(z))(d(z))^{1-\frac{1}{p}}|f(z)|^{p-1}(d(z))^{\frac{1}{p}}\|D_{f}(z)\|\\
\nonumber &&+p\sup_{z\in\mathbb{D}}(b(z))d(z)|f(z)|^{p}+p\sup_{z\in\mathbb{D}}(q(z))d(z)|f(z)|^{p-1}\\
&\leq&C_{8} (d(z))^{\frac{2}{p}}\|D_{f}(z)\|^{2}
+C_{9}(d(z))^{\frac{1}{p}}\|D_{f}(z)\|+C_{10}, \eeq where
$\displaystyle
C_{8}=2^{p-2}p^{2}\left(|f(0)|^{p-2}+C_{7}^{p-2}\right)$,
$\displaystyle
 C_{9}=2^{p-2}p\sup_{z\in\mathbb{D}}(a(z))\left(|f(0)|^{p-1}+C_{7}^{p-1}\right)
$ and $\displaystyle
C_{10}=2^{p-1}p\sup_{z\in\mathbb{D}}(b(z))\left(|f(0)|^{p}+C_{7}^{p}\right)+
2^{p-2}p\sup_{z\in\mathbb{D}}(q(z))\left(|f(0)|^{p-1}+C_{7}^{p-1}\right).
$ By the Cauchy-Schwarz inequality, we get \beq\label{eqt-1}
\left(\int_{\mathbb{D}}d^{\frac{1}{p}}(z)\|D(z)\|\,d\sigma(z)
\right)^{2}&\leq&\int_{\mathbb{D}}d^{\frac{2}{p}}(z)\|D(z)\|^{2}d\sigma(z)\int_{\mathbb{D}}
\,d\sigma(z)\\
\nonumber &=&\|f\|_{\mathcal{D}_{\gamma,2}}<\infty.
\eeq
Hence (\ref{eq-50}) and (\ref{eqt-1}) imply
\beq\label{eqr1}
\int_{\mathbb{D}}d(z)\Delta(|f(z)|^{p})\,d\sigma(z)&\leq&\int_{\mathbb{D}}\Big[C_{8}
(d(z))^{\frac{2}{p}}\|D_{f}(z)\|^{2}\\ \nonumber
&&+C_{9}(d(z))^{\frac{1}{p}}\|D_{f}(z)\|+C_{10}\Big]d\sigma(z)\\
\nonumber &\leq&C_{8} \|f\|_{\mathcal{D}_{\gamma,2}}+C_{9}
\|f\|_{\mathcal{D}_{\gamma/2,1}}+C_{10}\\ \nonumber
&<&\infty.
\eeq

%\bca \label{c-2} Let $p\in[2,4)$. \eca
In the case $p\in[2,4)$,  we let $F_{n}^{p}=(|f|^{2}+\frac{1}{n})^{p/2}$ for $n\in\{1,2,\ldots\}$.
We see that
$\Delta(F_{n}^{p})$ is integrable in $\mathbb{D}_{r}$. Then, by
(\ref{eq1.2x}), (\ref{eq-50}), (\ref{eqr1})  and Lebesgue's
dominated convergence theorem, we have
\beq\label{eq-t2} \nonumber
\lim_{n\rightarrow\infty}\int_{\mathbb{D}_{r}}d(z)\Delta(F_{n}^{p}(z))\,d\sigma(z)
&=&\int_{\mathbb{D}_{r}}d(z)
\lim_{n\rightarrow\infty}\big[\Delta(F_{n}^{p}(z))\big]\,d\sigma(z)\\
\nonumber &=&\frac{1}{2}\int_{\mathbb{D}_{r}}\Big[p(p-2)|
f(z)|^{p-4}|f(z)\overline{f_{z}(z)}+
f_{\overline{z}}(z)\overline{f(z)}|^{2}\\ \nonumber
&&+2p|f(z)|^{p-2}\big(|f_{z}(z)|^{2}+|f_{\overline{z}}(z)|^{2}\big)\\
\nonumber &&+p|f(z)|^{p-2}{\rm Re}\big(\overline{f(z)}\Delta
f(z)\big)\Big]d(z)\,d\sigma(z)\\
\nonumber &\leq&\int_{\mathbb{D}_{r}}\big[C_{8} d^{\frac{2}{p}}(z)\|D_{f}(z)\|^{2}\\
\nonumber
&&+C_{9}d^{\frac{1}{p}}(z)\|D_{f}(z)\|+C_{10}\big]\,d\sigma(z)\\
\nonumber
 &<&\infty. \eeq
Therefore, (\ref{eq-p1}) follows from the two cases.
% \ref{c-1} and \ref{c-2}.

%By , we conclude that for $p\in[2,\infty)$,

%\be\label{eq-o1}\int_{\mathbb{D}}d(z)\Delta(|f(z)|^{p})\,d\sigma(z)<\infty.\ee

Next we prove $f\in H^{p}_{g}(\mathbb{D}).$  As in the proof of Theorem~1.4 in \cite{CP}, for a fixed $r\in(0,1)$, since
$$\lim_{|z|\rightarrow r}\frac{\log r-\log |z|}{r-|z|}=\frac{1}{r},
$$
we see that there is an $r_{0}\in(0,r)$ satisfying
\be\label{eq-h1}
\log r-\log |z|\leq\frac{2}{r}(r-|z|)
\ee
for $r_{0}\leq|z|<r$.  Then it follows from
$\lim_{\rho\rightarrow0+}\rho\log (1/\rho )=0$
that
\beq\label{eq-x5}
\int_{\mathbb{D}_{r_{0}}}\Delta
(|f(z)|^{p})\log\frac{r}{|z|}\,d\sigma(z)&\leq&\int_{\mathbb{D}_{r_{0}}}\Delta
(|f(z)|^{p})\log\frac{1}{|z|}\,d\sigma(z)\\ \nonumber
&=&\int_{0}^{2\pi}\int_{0}^{r_{0}}\Delta (|f(\rho
e^{i\theta})|^{p})\rho\log\frac{1}{\rho}\,d\rho\, d\theta\\
\nonumber&<&\infty.
\eeq
Hence, by (\ref{eq1.2x}),  (\ref{eq-p1}), (\ref{eq-h1}) and (\ref{eq-x5}), we obtain
\begin{eqnarray*}
M_{p}^{p}(r,f)&=&|f(0)|^{p}+ \frac{1}{2}\int_{\mathbb{D}_{r}}\Delta
(|f(z)|^{p})\log\frac{r}{|z|}\,d\sigma(z)\\
&=&|f(0)|^{p}+\frac{1}{2}\int_{\mathbb{D}_{r_{0}}}\Delta
(|f(z)|^{p})\log\frac{r}{|z|}\,d\sigma(z)\\
&&+\frac{1}{2}\int_{\mathbb{D}_{r}\backslash
\mathbb{D}_{r_{0}}}\Delta
(|f(z)|^{p})\log\frac{r}{|z|}\,d\sigma(z)\\
 &\leq&|f(0)|^{p}+\frac{1}{2}\int_{\mathbb{D}_{r_{0}}}\Delta
(|f(z)|^{p})\log\frac{r}{|z|}\,d\sigma(z)\\
&&+ \int_{\mathbb{D}_{r}\backslash \mathbb{D}_{r_{0}}}\Delta
(|f(z)|^{p})\frac{(r-|z|)}{r}\,d\sigma(z)\\
&\leq&|f(0)|^{p}+\frac{1}{2}\int_{\mathbb{D}_{r_{0}}}\Delta
(|f(z)|^{p})\log\frac{1}{|z|}\,d\sigma(z)\\
&&+ \int_{\mathbb{D}\backslash \mathbb{D}_{r_{0}}}d(z)\Delta
(|f(z)|^{p})\,d\sigma(z)\\
 &<&\infty,
 \end{eqnarray*}
which implies that
%$$\sup_{0<r<1}M_{p}^{p}(r,f)<\infty.
%$$
$f\in H^{p}_{g}(\mathbb{D})$.  \qed

%Since Lemma \ref{lem-3} shows that the function $M_{p}^{p}(r,f)$ is
%increasing with respect to $r$ in $(0,1)$, we know that the limit
%$$\lim_{r\rightarrow1-}M_{p}(r,f)$$ does exist, which implies
%$f\in H^{p}$. The proof of the theorem is complete. \qed

\section{Lipschitz-type spaces  }\label{csw-sec3}
The following simple lemma is useful in the sequel.

\begin{lem}\label{lem-r1}
Let $\omega$ be a majorant and $\nu\in(0,1]$. Then for
$t\in(0,\infty)$, $\omega(\nu t)\geq\nu\omega(t)$.
\end{lem}
%\subsection*{Proof of Lemma \ref{lem-r1}}
\bpf Since $\omega(t)/t$ is decreasing on $t\in(0,\infty)$, we see that
$$\frac{\omega(\nu t)}{\nu t}\geq\frac{\omega( t)}{ t}
$$
and the desired conclusion follows.
\epf
%Then $\omega(\nu t)\geq\nu\omega(t)$. The proof of this lemma is complete.
%\qed

\subsection*{Proof of Theorem \ref{thm-1}} We first prove the sufficiency. For $r\in(0,1)$ and
$\theta\in[0,2\pi]$, let $w=z+re^{i\theta}$. Then
\begin{eqnarray*}
\|D_{f}(z)\|&=&\max_{\theta\in[0,2\pi]}\left|f_{x}(z)\cos\theta+f_{y}(z)\sin\theta\right|\\
&=&\max_{\theta\in[0,2\pi]}\left\{\lim_{r\rightarrow0+}\frac{|f(z+re^{i\theta})-f(z)|}{r}\right\}\\
&=&\max_{\theta\in[0,2\pi]}\left\{\lim_{r\rightarrow0+}\frac{|f(z)-f(w)|}{|z-w|}\right\}\\
%&\leq&\sup_{z,w\in\mathbb{D},z\neq
%w}\left\{\frac{|f(z)-f(w)|}{|z-w|}\right\}\\
&\leq&\lim_{r\rightarrow0+}\frac{C_{1}}{\omega\big(d^{s}(z)d^{\alpha-s}(z+re^{i\theta})}\\
&=&\frac{C_{2}}{\omega(d^{\alpha}(z))}.
\end{eqnarray*}

Next we  prove  the necessity. For $z, w\in\mathbb{D}$,  let
$\chi(t)=zt+(1-t)w,$ where $t\in[0,1].$ Since
\begin{eqnarray*}
1-|\chi(t)|\geq 1-t|z|-|w|+t|w| \geq (1-t)(1-|w|)=(1-t)d(w)
%\\
%&\geq&1-t+|w|(t-1)\\
%&=&(1-t)d(w)
\end{eqnarray*}
and similarly,  $1-|\chi(t)|\geq td(z)$,
%\begin{eqnarray*}
%1-|\chi(t)|&\geq&1-t|z|-|w|+t|w|\\
%&\geq&1-t|z|-(1-t)\\
%&=&td(z),
%\end{eqnarray*}
we see that
\be\label{eq-4}
(1-|\chi(t)|)^{\alpha-s}\geq(1-t)^{\alpha-s}d^{\alpha-s}(w)
\ee
and
\be\label{eq-5}
(1-|\chi(t)|)^{s}\geq t^{s}d^{s}(z).
\ee
By (\ref{eq-4}) and (\ref{eq-5}), we get
$$t^{s}(1-t)^{\alpha-s}d^{s}(z)d^{\alpha-s}(w)\leq(1-|\chi(t)|)^{\alpha},
$$
which implies
$$\omega\left(t^{s}(1-t)^{\alpha-s}d^{s}(z)d^{\alpha-s}(w)\right)\leq\omega\left((1-|\chi(t)|)^{\alpha}\right)
=\omega\big(d^{\alpha}(\chi(t))\big).
$$
Hence, for $z, w\in\mathbb{D}$ with $z\neq w$, by Lemma
\ref{lem-r1}, we know that there is a positive constant $C$ such
that
\begin{eqnarray*}
|f(z)-f(w)|&=&\left|\int_{0}^{1}\frac{df}{dt}(\chi(t))\,dt\right|\quad (\zeta =\chi(t))\\
&=&\left|(z-w)\int_{0}^{1}f_{\zeta}(\chi(t))\,dt+(\overline{z}-\overline{w})\int_{0}^{1}f_{\overline{\zeta}}(\chi(t))\,dt\right|\\
&\leq&|z-w|\int_{0}^{1}\|D_{f}(\chi(t))\|\,dt\\
&=&|z-w|\int_{0}^{1}\frac{\|D_{f}(\chi(t))\|}{\omega\left(d^{\alpha}(\chi(t))\right)}
\omega\left(d^{\alpha}(\chi(t))\right)\,dt\\
&\leq&C|z-w|\int_{0}^{1}\frac{dt}{\omega\left(d^{\alpha}(\chi(t))\right)}
\\
%\end{eqnarray*}
%\begin{eqnarray*}
&\leq&C|z-w|\int_{0}^{1}\frac{dt}{\omega\left(t^{s}(1-t)^{\alpha-s}d^{s}(z)d^{\alpha-s}(w)\right)}
%\\
\end{eqnarray*}
\begin{eqnarray*}
&\leq&\frac{C|z-w|}{\omega\left(d^{s}(z)d^{\alpha-s}(w)\right)}\int_{0}^{1}\frac{dt}{(1-t)^{\alpha-s}t^{s}}\\
&=&\frac{C|z-w|}{\omega\left(d^{s}(z)d^{\alpha-s}(w)\right)}\mathbf{B}(1-s,1+s-\alpha),
\end{eqnarray*}
where $\mathbf{B}(\cdot,\cdot)$ denotes the Beta function. Thus,
there is a positive constant $C_{1}=C\mathbf{B}(1-s,1+s-\alpha)$ such that for all $z$
and $ w$ with $z\neq w$,
$$\frac{|f(z)-f(w)|}{|z-w|}\leq \frac{C_{1}}{\omega\big(d^{s}(z)d^{\alpha-s}(w)\big)}.
$$
%Since
%$0\leq s<1$ and $s\leq\alpha<s+1$, we know that
%$\mathbf{B}(1-s,1+s-\alpha)<\infty.$ Hence there is a positive
%constant $C_{1}=C\mathbf{B}(1-s,1+s-\alpha)$ such that for all $z$
%and $ w$ with $z\neq w$,
%$$\frac{|f(z)-f(w)|}{|z-w|}\leq \frac{C_{1}}{\omega\big(d^{s}(z)d^{\alpha-s}(w)\big)}.
%$$
The proof of this theorem is complete. \qed

\begin{Lem}$($\cite[Lemma 2.2]{CPMW}$)$\label{A}
Suppose that $f$ is a harmonic mapping in
$\overline{\mathbb{D}}(a,r)$, where $a\in\mathbb{C}$ and $r>0$. Then
$$\|D_{f}(a)\|\leq\frac{2}{\pi r}\int_{0}^{2\pi}|f(a+re^{i\theta})-f(a)|\,d\theta.
$$
\end{Lem}

\subsection*{Proof of Theorem \ref{thm-2}}
We first prove the sufficiency. By Lemma \Ref{A}, for
$\rho\in(0,d(z)]$,
$$\|D_{f}(z)\|\leq\frac{2}{\pi \rho}\int_{0}^{2\pi}\big|f(z+\rho e^{i\theta})-f(z)\big|\,d\theta,
$$
which gives
$$\int_{0}^{r}\rho^{2}\|D_{f}(z)\|d\rho\leq\frac{2}{\pi }\int_{0}^{r}\left(\rho\int_{0}^{2\pi}
|f(z+\rho e^{i\theta})-f(z)|d\theta\right)d\rho,
$$
where $r=d(z).$ Then
\begin{eqnarray*}
\|D_{f}(z)\|&\leq&\frac{6}{\pi
r^{3}}\int_{\mathbb{D}(z,r)}|f(z)-f(\zeta)|\,dA(\zeta)\\
&=&\frac{6}{r |\mathbb{D}(z,r)|}\int_{\mathbb{D}(z,r)}|f(z)-f(\zeta)|\,dA(\zeta)\\
%&\leq&\frac{6C_{2}r^{\alpha-1}}{\omega(r^{\alpha})}\\
&\leq&\frac{6C_{2}}{\omega(r^{\alpha})}.
\end{eqnarray*}

Now we prove the necessity. Since
$f\in\mathcal{L}_{\infty,\omega}\mathcal{B}^{0}_{\alpha}(\mathbb{D})$,
we see that there is a positive constant $C$ such that
\be\label{eq-11}
\|D_{f}(z)\|\leq\frac{C}{\omega(d^{\alpha}(z))}.
\ee
For $z,w\in\mathbb{D}$ and $t\in[0,1]$, if $d(z)>t|z-w|$, then, by (\ref{eq-11}), we get
\begin{eqnarray*}
|f(z)-f(w)|&\leq&|z-w|\int_{0}^{1}\|D_{f}(z+t(w-z))\|\,dt\\
&\leq&C|z-w|\int_{0}^{1}\frac{dt}{\omega\big(d^{\alpha}(z+t(w-z))\big)}\\
&\leq&C|z-w|\int_{0}^{1}\frac{dt}{\omega\left(\big(d(z)-t|z-w|\big)^{\alpha}\right)}\\
&=&C\int_{0}^{|z-w|}\frac{dt}{\omega\left(\big(d(z)-t\big)^{\alpha}\right)},
\end{eqnarray*}
which implies
\begin{eqnarray*}
\frac{1}{
|\mathbb{D}(z,r)|}\int_{\mathbb{D}(z,r)}|f(z)-f(\zeta)|\,dA(\zeta)
&\leq&\frac{C}{|\mathbb{D}_{r}|}\int_{\mathbb{D}_{r}}\left(\int_{0}^{|\xi|}
\frac{dt}{\omega\left(\big(d(z)-t\big)^{\alpha}\right)}\right)dA(\xi)
\\
&=&\frac{2C}{r^{2}}\int_{0}^{r}\rho\left(\int_{0}^{\rho}
\frac{dt}{\omega\left(\big(d(z)-t\big)^{\alpha}\right)}\right)d\rho
%\\
\end{eqnarray*}
\begin{eqnarray*}
&\leq&\frac{2C}{r^{2}}\int_{0}^{r}\left(\int_{t}^{r}\rho
d\rho\right)\frac{dt}{\omega\left(\big(r-t\big)^{\alpha}\right)}\\
&=&\frac{2C}{r}\int_{0}^{r}\frac{\big(r-t\big)^{\alpha}}{\omega\left(\big(r-t\big)^{\alpha}\right)}
\big(r-t\big)^{1-\alpha}\,dt\\
&\leq&\frac{2Cr^{\alpha-1}}{\omega(r^{\alpha})}\int_{0}^{r}\big(r-t\big)^{1-\alpha}\,dt\\
&=&C_{2}\frac{r}{\omega(r^{\alpha})},
\end{eqnarray*}
where $C_{2}=\frac{2C}{2-\alpha}.$ The proof of this theorem is
complete. \qed

\subsection*{Proof of Theorem \ref{thm-3}}
We first prove the necessity. For $r\in(0,1),$ let $F(z)=f(rz).$ By
the proof of necessity part of Theorem \ref{thm-1}, we see that
there is a positive constant $C$ such that
\be\label{eq-13}
\frac{\big|(F(z)-f(z))-(F(w)-f(w))\big|\omega\left(d^{s}(z)d^{\alpha-s}(w)\right)}{|z-w|}\leq
C\|f-F\|_{\mathcal{L}_{\infty,\omega}\mathcal{B}^{0}_{\alpha}(\mathbb{D})}.
\ee
Since $\omega(t)/t$ is non-increasing for $t>0$, we know that there is a positive constant $C$ such that
\begin{eqnarray*}
\frac{\big|F(z)-F(w)\big|\omega\left(d^{s}(z)d^{\alpha-s}(w)\right)}{|z-w|}&=&
\frac{r\big
|F(z)-F(w)\big|\omega\left(d^{s}(rz)d^{\alpha-s}(rw)\right)}{|rz-rw|}\\
&&\times\frac{\omega\left(d^{s}(z)d^{\alpha-s}(w)\right)}{\omega\left(d^{s}(rz)d^{\alpha-s}(rw)\right)}\\
&\leq&Cr\|f\|_{\mathcal{L}_{\infty,\omega}\mathcal{B}^{0}_{\alpha}(\mathbb{D})}
\frac{\omega\left(d^{s}(z)d^{\alpha-s}(w)\right)}{\omega\left(d^{s}(rz)d^{\alpha-s}(rw)\right)}\\
&=&Cr\|f\|_{\mathcal{L}_{\infty,\omega}\mathcal{B}^{0}_{\alpha}(\mathbb{D})}
\frac{\frac{\omega\left(d^{s}(z)d^{\alpha-s}(w)\right)}{d^{s}(z)d^{\alpha-s}(w)}}
{\frac{\omega\left(d^{s}(rz)d^{\alpha-s}(rw)\right)}{d^{s}(rz)d^{\alpha-s}(rw)}}\\
&&\times\frac{d^{s}(z)d^{\alpha-s}(w)}{d^{s}(rz)d^{\alpha-s}(rw)}
\\
&\leq&
%Cr\|f\|_{\mathcal{L}_{\infty,\omega}\mathcal{B}^{0}_{\alpha}(\mathbb{D})}\frac{d^{s}(z)d^{\alpha-s}(w)}{d^{s}(rz)d^{\alpha-s}(rw)} \\
%\end{eqnarray*}
%\begin{eqnarray*}
%&=&
Cr\|f\|_{\mathcal{L}_{\infty,\omega}\mathcal{B}^{0}_{\alpha}(\mathbb{D})}\left(\frac{d(z)}{d(rz)}\right)^{s}
\left(\frac{d(w)}{d(rw)}\right)^{\alpha-s}\\
&\leq&Cr\|f\|_{\mathcal{L}_{\infty,\omega}\mathcal{B}^{0}_{\alpha}(\mathbb{D})}\left(\frac{d(z)}{d(rz)}\right)^{s}.
\end{eqnarray*}
By using the triangle inequality, we have
$$\sup_{ z\neq
w}\left\{\frac{|f(z)-f(w)|\omega\big(d^{s}(z)d^{\alpha-s}(w)\big)}{|z-w|}\right\}\leq
C\|f-F\|_{\mathcal{L}_{\infty,\omega}\mathcal{B}^{0}_{\alpha}(\mathbb{D})}+
Cr\|f\|_{\mathcal{L}_{\infty,\omega}\mathcal{B}^{0}_{\alpha}(\mathbb{D})}\left(\frac{d(z)}{d(rz)}\right)^{s}.
$$
In the above inequality, first letting $|z|\rightarrow1-$  and then
letting $r\rightarrow1-$, we get the desired result.

Next we begin to prove the sufficiency. Suppose (\ref{eq12}) holds.
For all $\epsilon>0$, there is a $\delta\in(0,1)$ such that
$$\sup_{ w\in\mathbb{D},z\neq
w}\left\{\frac{|f(z)-f(w)|\omega\big(d^{s}(z)d^{\alpha-s}(w)\big)}{|z-w|}\right\}<\epsilon,
$$
whenever $|z|>\delta$. Let $w$ tend to $z$ in the radial direction,
we obtain
$$\|D_{f}(z)\|\omega\big(d^{\alpha}(z)\big)\leq\epsilon
$$
whenever $|z|>\delta$, which yields
$f\in\mathcal{L}_{\infty,\omega}\mathcal{B}^{0}_{\alpha}(\mathbb{D})$.
  \qed

\section{Composition operators}\label{csw-sec4}

Given $f\in \mathcal{A}(\mathbb{D})$, the Littlewood-Paley $g$-function is
defined as follows
$$g(f)(\zeta)=\left(\int_{0}^{1}|f'(r\zeta)|^{2}(1-r)\,dr\right)^{\frac{1}{2}},~\zeta\in\partial\mathbb{D}.
$$
By \cite[Theorems 3.5 and 3.19]{Zy}, we know that $f\in
H^{p}(\mathbb{D})$ if and only if $g(f)\in H_{g}^{p}(\mathbb{D})$
for $p>1.$

\subsection*{Proof of Theorem \ref{thm-8}} We first prove that
(1)$\Longrightarrow$(2). Applying \cite[Lemma 1]{AD} and Lemma
\ref{lem-5}, we see that there are two functions
$f_{1},~f_{2}\in\mathcal{L}\mathcal{B}^{\beta}_{\alpha}(\mathbb{D})$
such that for $z\in\mathbb{D},$
\be\label{eq-j1}
|f_{1}'(z)|^{2}+|f_{2}'(z)|^{2}\geq
d^{-2\alpha}(z)\left(\log\frac{e}{d(z)}\right)^{-2\beta}.
\ee
Since for $k=1,2$, $C_{\phi}(f_{k})\in H^{2}(\mathbb{D})$, by
(\ref{eq-j1}), we conclude that
\begin{eqnarray*}
\infty&>&\|g(C_{\phi}(f_{1}))\|^{2}_{2}+\|g(C_{\phi}(f_{2}))\|^{2}_{2}\\
&=&\frac{1}{2\pi}\int_{0}^{2\pi}\int_{0}^{1}\left(|f_{1}'(\phi(r\zeta))|^{2}+|f_{2}'(\phi(r\zeta))|^{2}\right)
|\phi'(r\zeta)|^{2}(1-r)\,dr\,d\theta\\
&\geq&\frac{1}{2\pi}\int_{0}^{2\pi}\int_{0}^{1}
\frac{|\phi'(re^{i\theta})|^{2}}{d^{2\alpha}(\phi(re^{i\theta}))}
\left(\log\frac{e}{d(\phi(re^{i\theta}))}\right)^{-2\beta}(1-r)\,dr\,d\theta,
\end{eqnarray*}
which shows that (1)$\Longrightarrow$(2).

Next we  prove (2)$\Longrightarrow$(1). For
$f\in\mathcal{L}\mathcal{B}^{\beta}_{\alpha}(\mathbb{D})$ and
$\zeta\in\partial\mathbb{D}$, we get
\begin{eqnarray*}
g^{2}(C_{\phi}(f))(\zeta)&=&\int_{0}^{1}|\big(C_{\phi}(f)(r\zeta)\big)'|^{2}(1-r)\,dr\\
&=&\int_{0}^{1}|f'(\phi(r\zeta))|^{2} |\phi'(r\zeta)|^{2}(1-r)\,dr
\\
%\end{eqnarray*}
%\begin{eqnarray*}
&=&\int_{0}^{1}|f'(\phi(r\zeta))|^{2}d^{2\alpha}(\phi(re^{i\theta}))\left(\log\frac{e}{d(\phi(re^{i\theta}))}\right)^{2\beta}\\
&&\times|\phi'(r\zeta)|^{2}d^{-2\alpha}(\phi(re^{i\theta}))\left(\log\frac{e}{d(\phi(re^{i\theta}))}\right)^{-2\beta}(1-r)\,dr\\
&\leq&\|f\|_{\mathcal{L}\mathcal{B}^{\beta}_{\alpha}(\mathbb{D})}^{2}\int_{0}^{2\pi}\int_{0}^{1}
\frac{|\phi'(re^{i\theta})|^{2}}{d^{2\alpha}(\phi(re^{i\theta}))}
\left(\log\frac{e}{d(\phi(re^{i\theta}))}\right)^{-2\beta}(1-r)\,dr,
\end{eqnarray*}
which yields (2)$\Longrightarrow$(1), whence $g(C_{\phi}(f))\in
H_{g}^{2}(\mathbb{D})$. The proof of this theorem is complete. \qed

%\bigskip

%{\bf Acknowledgements:} This research was  partly  supported by NSF
%of China (No. 11071063). This work was also supported in part by the
%Construct Program of the Key Discipline in Hunan Province (No.
%[2011] 76) and the Start Project of Hengyang Normal University (No.
%12B34).
\subsection*{Acknowledgments} This research was partly supported by NSF
of China (No. 11326081), the Construct Program of the Key Discipline in Hunan Province,
Academy of Finland (No. 269260) and National Natural Science Foundation of China.

\normalsize

\end{document}